\documentclass[]{amsart}

\usepackage{amsfonts}

\newtheorem{theorem}{Theorem}[section]        
\newtheorem{proposition}[theorem]{Proposition} 

\newtheorem{lemma}[theorem]{Lemma}

\theoremstyle{definition}

\newtheorem{definition}[theorem]{Definition}

\theoremstyle{remark}
\newtheorem{remark}[theorem]{Remark}

\numberwithin{equation}{section}

\newcommand{\rn} {{\mathbb R}^n}
\newcommand{\rstar} {{\mathbb R}_*}

\newcommand{\tr} {\textup{Tr}}
\newcommand{\bbar}{\left( \begin{array}}
\newcommand{\ebar}{\end{array} \right) }
\newcommand{\leftb}{\left< \,}
\newcommand{\rightb}{\, \right> }
\newcommand{\bdm}{\begin{displaymath}}
\newcommand{\edm}{\end{displaymath}}
\newcommand{\beq}{\begin{equation}}
\newcommand{\beqa}{\begin{eqnarray}}
\newcommand{\beqas}{\begin{eqnarray*}}
\newcommand{\eeq}{\end{equation}}
\newcommand{\eeqa}{\end{eqnarray}}
\newcommand{\eeqas}{\end{eqnarray*}}
\newcommand{\dd}{\textup{d}}
\newcommand{\pip}{\pi_{\mathcal{P}}}
\newcommand{\pipa}{\pi_{\mathcal{P}_a}}
\newcommand{\pin}{\pi_{\mathcal{N}}}
\newcommand{\gtilde}{\widetilde{\mathcal{G}}}
\newcommand{\gtildesig}{\widetilde{\mathcal{G}}_{\sigma}}

\newcommand{\gtildesigr}{{\mathbb R}\gtildesig}
\newcommand{\gltilde}{\widetilde{\textup{gl}}(m,{\mathbb C})}
\newcommand{\pa}{\mathcal{P}_a}
\newcommand{\na}{\mathcal{N}_a}

\title[$k$-symmetric AKS systems and flat 
immersions in spheres]
{$k$-symmetric AKS systems and Flat Immersions into Spheres}

\author{David Brander}
\address{Department of Mathematics\\ Faculty of Science\\Kobe University\\1-1, Rokkodai, Nada-ku, Kobe 657-8501\\ Japan}
\email{brander@math.kobe-u.ac.jp}


\begin{document}

\begin{abstract}
We define a large class of integrable nonlinear PDE's, \emph{$k$-symmetric
AKS systems}, whose solutions evolve on finite dimensional subalgebras
of loop algebras, and linearize on an associated algebraic curve.
We prove that periodicity of the associated algebraic data implies a
type of quasiperiodicity for the solution, and show that the problem of
isometrically immersing $n$-dimensional Euclidean space into a sphere
of dimension $2n-1$  can be addressed via this scheme, producing 
infinitely many real analytic solutions.
\end{abstract}

\maketitle

\section{Introduction}
In this paper we apply the Adler-Kostant-Symes (AKS) theory of integrable
systems to the problem of isometrically immersing Euclidean $n$-space, $E^n$,
into a sphere of dimension $2n-1$.  We show how the theory can be used
to produce families of complete solutions which linearize on an associated
algebraic variety, and study the implications of periodicity of the 
associated solutions. It is the author's hope that this example will be
of value to differential geometers wishing to understand how the AKS
theory is applied to problems in geometry, as it demonstrates 
many of the main features of the theory in a particularly simple and
uncluttered setting.  It should be mentioned here that there are other
methods from the theory of loop groups, such as dressing of solutions,
which could be applied to our problem, but are not covered here.

The codimension of $n-1$ for the immersion is of interest 
because it  is critical in
the sense that there is no local isometric embedding of $E^n$ into
$S^k$ for $k< 2n-1$ (see \cite{spivak5}, p. 195), whilst for
$k = 2n-1$ we at least
have the Clifford tori,  products $S^1 \times S^1 ... \times S^1$,
of $n$ appropriately scaled circles. In fact the theory described in
this paper also applies to codimension greater than  $n-1$, if we
add the condition that the normal bundle is flat.

The study of isometric immersions from $E^n$ into $S^{2n-1}$ has a
history going back at least as far as 1896, when Bianchi \cite{bianchi}
 studied flat surfaces in the 3-sphere. These can be classified roughly 
using  the existence of special ``asymptotic coordinates'' as well as
 the group structure on $S^3$ (see \cite{spivak4}, pp. 139-163).
More recently J. Weiner \cite{weiner} and Y. Kitagawa \cite{kitagawa} were
able to extend this to classify solutions with two periods, that is 2-tori
in $S^3$. See also \cite{dadoksha}.

When $n>2$ however, one cannot use the same methods: although 
one does have asymptotic coordinates \cite{moore}, there is no group
structure on a higher dimensional sphere.  The system of PDE associated
with the integrability condition, the so-called `generalized wave
equation' has been studied to some extent by Tenenblat \cite{tenenblat}
and others, though not from the AKS viewpoint.  
 
On the other hand, over the past two decades, methods based on
the work of Adler and Van Moerbeke
 \cite{adlervanmoerbeke1}, \cite{adlervanmoerbeke2},
Kostant \cite{kostant1979} and Symes \cite{symes1980} have been
applied successfully to problems of producing and sometimes
classifying certain special submanifolds.
 Following work by Hitchin \cite{hitchin1986},  
Pinkall and Sterling \cite{pinkallsterling} gave a classification 
of all constant mean curvature 2-tori in $E^3$ using the so-called
``finite gap'' method. 
 The finite gap method was also used to complete a 
classification of minimal 2-tori in $S^4$ in \cite{fpps}. All this,
and other work, such as Uhlenbeck's study
 \cite{uhlenbeck1989} of harmonic maps into Lie groups via loop groups,
led to a general scheme for obtaining harmonic tori in symmetric spaces,
described in \cite{bfpp}. Recently there has been
work by Terng and others
on a programme to identify the submanifolds associated to a large
class of related integrable systems  \cite{terng1997}, \cite{terng2002a}.

Closely related to our target problem,  Ferus and
Pedit, \cite{feruspedit1996}, \cite{feruspedit1996II},  
modified the AKS theory to obtain  large families of local isometric
immersions of a space form $M_c^n$ into another $M_{\tilde{c}}^{2n-1}$
(the subscript is the constant sectional curvature of the manifold),
where $c$ and $\tilde{c}$ are both non-zero. These were special cases
of so-called \emph{curved flats}, which were later also studied in
 \cite{ternguhlenbeck2000}, where dressing actions
were used to produce families of solutions. 

Below we show that the case $c =0$ can be studied in a similar manner.
We begin by introducing some known theory of integrable systems in a 
form which is intended to be easily applied to problems in differential
geometry. Specifically, we  define a large class of 
integrable systems, which we shall call 
\emph{$k$-symmetric AKS systems}.  These systems
are defined by a Lax equation 
\beq  \label{lax1}
\dd X = [X,A],
\eeq
 where $X(t)$,
for $t$ in $\rn$, takes its
values in a loop algebra of matrix-valued polynomials in some auxiliary
parameter $z$. They have the 
property that they are of \emph{finite type}, in the sense that if the
initial condition 
$X(0)$ is a polynomial, then $X(t)$ is also a polynomial of the
same degree for all $t$.  

There is 
associated to $X(t)$ a line bundle $L(t)$ over a (constant) algebraic curve
$\Sigma$.  The evolution of $L(t)$ in the Picard group of line bundles
(of fixed degree) over $\Sigma$ is linear, which, in a sense, 
linearizes the Lax equation. 

For applications in submanifold theory, we are actually interested
in a map $F$ into a Lie group associated to the loop algebra mentioned
above. $F$  can be characterized here by the equation
\beqas
A = F^{-1} \dd F, \\
F(0) = I.
\eeqas
We will call $F$ a \emph{$k$-symmetric map} and
 $X$ a \emph{Killing field for $F$}.

Now the question of periodicity is easily addressed for the
linearized solution $L(t)$, so we would like to know what this
means for the function $F(t)$.
In this paper we prove, for a large subclass of $k$-symmetric AKS systems,
called \emph{simple}, the following:
\begin{proposition} \label{proposition1}
Suppose there exists an element $P$ in $\rn$ such that the
 line bundle $L(t)$ associated to a solution $X(t)$ of (\ref{lax1})
satisfies 
\bdm
L(t+P) = L(t),
\edm 
for all $t$ in $\rn$.  Then there exists a constant 
matrix $B \in \textup{GL}(2n, {\mathbb C})$, such that the following 
two equations hold for all $t$ in $\rn$:

\beqa
X(t+P) = B^{-1}X(t) \, B, \\
F(t+P) = F(P) \, B^{-1}F(t) \, B.
\eeqa

\end{proposition}

In Section \ref{applications} we study the problem 
of producing  isometric immersions
from $E^n$ into $S^{2n-1}$ and show that it is a 2-symmetric 
AKS system.  Specifically, we prove Theorem \ref{theorem1} below.

To state the result, we
first need some definitions. If $f: \rn \to S^{2n-1}$ is a smooth
map, then we will say $F: \rn \to SO(2n)$ is an \emph{adapted frame}
for $f$ if the $(n+1)$'th column of $F$ is $f$ and if the derivative
$\dd f$ has no component in the directions given by the last $n$
columns.  This means that if $f$ is an immersion then the first
$n$ columns span the tangent space and the last $n$ the normal space
to the immersion, considered as a map into ${\mathbb R }^{2n}.$ 

Let 
$so(2n,{\mathbb C})$ denote the Lie algebra of complex-valued skew
symmetric $2n \times 2n$ matrices. Consider the
 involution $\sigma = Ad(Q)$ on $so(2n, {\mathbb C})$, where 
$Q$ is the matrix  $\textup{diag}(
      I_{n \times n}, -I_{n \times n})$,
and $I_{n \times n}$ is the $n \times n$ identity matrix.
Let $\mathcal{V}_0 \subset so(2n,{\mathbb C})$ be the fixed point set
 of $\sigma$ and $\mathcal{V}_1$ the $-1$ eigenspace.
Consider the Lie algebra of matrix valued Laurent polynomials
\bdm
\gtildesig := \{ \sum_{i=\alpha}^{\beta} X_i z ^i ~|~ 
X_i \in \mathcal{V}_{i ~ \textup{(mod2)}} \},
\edm
and define for each $d \geq 0$ the real vector sub-space
\bdm 
(\gtildesigr )_d^1 := \{ \sum_{i = -d}^1 X_i z ^i ~|~ 
X_i \in   so(2n,{\mathbb R}) \cap \mathcal{V}_{i ~\textup{(mod2)}} \}.
\edm

\begin{theorem} \label{theorem1}
For every initial condition $X(0)$ in $(\gtildesigr )_d^1$ we have:
\begin{enumerate}
\item
A family $f^z$, parameterized by $z$ in the non-zero real numbers
$\rstar$, of complete real-analytic maps from $\rn$ into $S^{2n-1}$, together
with a family of adapted frames $F^z$ for $f^z$, which comes from
a 2-symmetric AKS system. For a generic initial condition, $f^z$ is
an immersion in a neighbourhood of $0$, and 
 wherever $f^z$ is immersive, the induced
curvature on its image is zero. 
\item
Similar families $f^z_P$ and $F^z_P$, which have the additional property
that both the tangent and normal frames given by $F^z_P$ are parallel.
\end{enumerate}
\end{theorem}

\begin{remark}
 One obtains a similar result, without the completeness, for 
 flat immersions into hyperbolic space, by substituting the group
$SO_{-1}(2n)$ for $SO(2n)$ (c.f. \cite{feruspedit1996}).
\end{remark}

\section{Lax Equations and Loop Groups} \label{section2}

In this section we outline the AKS theory as it applies to
$k$-symmetric AKS systems.  The main ideas here can be traced
to the papers \cite{adlervanmoerbeke1} and \cite{adlervanmoerbeke2} and
their predecessors.

\subsection{Algebraically Completely Integrable Systems}
First we sketch the most important facts about Lax equations and their
spectral data, more or less following \cite{fpps} and 
\cite{hitchin1999} for our definitions.

\subsubsection{Lax Equations} \label{subsectionLax}
A \emph{Lax equation} is a differential equation which can be 
written in the form:
\bdm
 \textup{d} X = [X,A], 
\edm
where $A = \sum_{i=1}^n A_i \dd t_i$, and the square brackets represents
the commutator $XA - AX$. 
Here we assume $X$ and $A_i$ are
complex valued matrices.

The integrability condition which ensures the existence of a local
solution, \emph{Frobenius integrability}, 
is $\dd ^2 X =0$ and  will be addressed below.

 Lax equations are sometimes called \emph{isospectral} because the
spectrum of a solution $X$,
\bdm
 \{ h \in {\mathbb C} | ~ \det (h I - X) = 0 \} ,
\edm
is \emph{constant} with respect to $t$.  This well-known fact can be verified
by observing that the coefficients of the characteristic  equation of
$X$ can all be written in terms of traces of powers of $X$ and then 
checking that $\dd(\tr X^k) = k \tr (X^{k-1} [X, A] = 0$.


\subsubsection{The Spectral Curve} \label{subsectspectcurve}
Let us now allow $X$ and $A$ to depend meromorphically on a
complex parameter $z$:
\beq  \label{lax}
 \textup{d} X(z,t) = [X(z,t), A(z,t)], ~~~~~~ z \in {\mathbb C}, ~ t \in {\mathbb R}^n.
\eeq
The case we will consider is that $X$ and $A_i$ are Laurent
polynomials in $z$, that is
\bdm
X(t),~A_i(t) \, \in \, \gltilde := \textup{gl}(m,{\mathbb C}) [z, z^{-1}].
\edm

  Rather than a discrete set, the spectrum of $X$ is now an
algebraic curve, called
the \emph{spectral curve} $\Sigma$ of $X$, defined as the algebraic completion of
the set: 
\bdm
\Sigma_0 = \{ (z, w) \in {\mathbb C}^2 |~  \det (w I - X(z)) = 0 \}.
\edm
$X$ is called a \emph{regular element} of \, $\gltilde$ if: 
\begin{enumerate}
\item
$\Sigma$ is a smooth algebraic curve.
\item
The projection $z: \Sigma \to {\mathbb CP}^1$, given by $[(z,w) \mapsto z]$,
is an $m$-fold branched covering with
simple branch points only and no branch points over the points $0$
and $\infty$.
\end{enumerate}
 Note that the last requirement implies that if $X = \sum_a^b X_i z^i$ 
then $X_a$ and $X_b$ both have
$m$ distinct eigenvalues. See \cite{fpps} for more details.

 We will always 
assume that $X(z,0)$ (and therefore $X(z,t)$) is regular, and this
is true for a generic element of $\gltilde$.

\subsubsection{The Eigenspace Bundle} \label{subsecteigenspacebundle}
Unlike the eigenvalues, the eigenvectors of
 $X(z,t)$ are not constant with respect to $t$.
They determine a family of line bundles
$L(t)$ over $\Sigma$, where $L(t)$ is the sub-bundle of the trivial
bundle $\Sigma \times {\mathbb C}^m$ whose fibre over $(z,w)$ is the
eigenspace corresponding to the eigenvalue $w$.

\subsubsection{Algebraic Complete Integrability}

 Considering that a matrix is determined by its eigenvalues and
eigenvectors, it is reasonable to study the solution $X(z,t)$ in terms of 
its \emph{spectral data} $\{ \Sigma, L(t), z \}$. 

 The degree $d$ of the line bundle $L(t)$, being an integer,
is constant by continuity, so $L(t)$ takes its values in
$\textup{Pic}^d (\Sigma)$, the space of equivalence classes of 
holomorphic line bundles over $\Sigma$ of degree
$d$.  It is well-known (see \cite{hitchin1999})
that this is a complex torus 
$\frac{{\mathbb C}^g}{{\mathbb Z}^{2g}}$, where $g$ is the genus of $\Sigma$,
and the Lax equation (\ref{lax}) is called \emph{algebraically completely
integrable} if $L(t)$ follows a straight line motion in this torus.
A general study of the question of algebraic 
complete integrability for (1-dimensional)
Lax equations can be found in \cite{griffiths}.

\subsubsection{Reconstructing $X(t)$ from its Spectral Data} \label{reconstruct}
Suppose given spectral data  $\{ \Sigma, L(t), z \}$ for some 
(unspecified) regular $X(t)$.
 The line bundle $L(t) \in \textup{Pic}^d (\Sigma)$, is only known up 
to holomorphic  equivalence. In other
words, the way $L(t)$ sits as a sub-bundle of $\Sigma \times {\mathbb C}^m$
is not a priori given to us. We need to know this in order to 
 recover the  solution
$X(t)$. The following result from \cite{fpps} is essentially what allows
one to do this:
\begin{proposition}  \label{dimprop}
Let $X \in \gltilde$ be regular and $L \to \Sigma$ its eigenspace
bundle with dual bundle $L^*$, and $\Gamma (L^*)$ the space of
holomorphic sections of $L^*$. Then 
\bdm
\textup{dim} \, \Gamma (L^*) = m.
\edm
\end{proposition}
This allows us to identify $(\Gamma (L^*))^*$ with ${\mathbb C}^m$ and thus
think of $L$ as a sub-bundle of $\Sigma \times {\mathbb C}^m$.  The identification
is not unique, so what we are able to reconstruct in general is of the
form:
\bdm
Y(z,t) = G^{-1}(t) X(z,t) G(t),
\edm
where $G(t)$ is an an element of $\textup{GL}(m,{\mathbb C})$, and is
independent of $z$.  
It is therefore reasonable to expect that some kind of normalization
condition on $X(t)$ may be needed in order to reconstruct it uniquely
from the spectral data.


\section{$k$-symmetric AKS systems}  \label{systemsoffinitetype}

Consider the Lie algebra 
$\widetilde{\textup{gl}}(m,{\mathbb C}) := \textup{gl}(m,{\mathbb C})[z,z^{-1}]$ of 
Laurent polynomials in $z$, which is given
the Lie bracket
\bdm
[\sum _i X_i z^i, \sum _j Y_j z^j] := \sum _{i,j} [X_i,Y_j] z^{i+j}.
\edm
$k$-symmetric AKS systems are maps into certain subalgebras 
$\gtilde_{\sigma}$ of $\widetilde{\textup{gl}}(m,{\mathbb C})$ which
arise as follows: suppose $\mathcal{G}$ 
is a Lie subalgebra of $\textup{gl}(m,{\mathbb{C}})$, and 
$\sigma$ is an order $k$ automorphism of $\mathcal{G}$, that is, 
$\sigma : \mathcal{G} \to \mathcal{G}$ is linear and satisfies 
$\sigma \, [X,Y] = [\sigma \, X, \sigma \, Y ]$ and $\sigma ^k = \textup{id}$.
 Let $\mathcal{V}_i$ be
the eigenspace of $\sigma$ corresponding to the eigenvalue 
$\zeta^i$, where $\zeta$ is a primitive $k$'th root of unity. 
If $X_i \in \mathcal{V}_i$, then
\beqas
\sigma [X_i, X_j] &=& [\sigma X_i, \, \sigma X_j] \\
&=& (\zeta)^{i+j} [X_i,X_j],
\eeqas
which implies the commutation relations 
\bdm
[\mathcal{V}_i, \mathcal{V}_j ] \subset \mathcal{V}_{i+j \, (\textup{mod} \,k)}.
\edm
It follows that the set
\bdm
\gtilde_{\sigma} := \{ X \in \mathcal{G} [z, z^{-1}] ~|~ X_i \in 
\mathcal{V}_{i ~ (\textup{mod} \,k)} \}
\edm
is closed under the Lie bracket, and therefore a Lie subalgebra of 
$\widetilde{\textup{gl}}(m,{\mathbb C})$.

 Now let $\gtilde$ be any subalgebra of $\gltilde$ and suppose we have
a decomposition
\beq \label{decomposition}
\gtilde = \mathcal{P} \oplus \mathcal{N}
\eeq
into two subalgebras which consist
of polynomials in $z$ and $z^{-1}$ respectively.  
An example of such a decomposition, which we will call \emph{simple}, is 
\beqa \label{simpledecomp}
\mathcal{P}  := \{ X \in \gtilde ~|~ X_i = 0,~ i \leq 0 \}, \\
\mathcal{N}  := \{ X \in \gtilde ~|~ X_i = 0,~ i > 0 \}. \nonumber
\eeqa
Other examples are obtained similarly by stipulating that elements
of $\mathcal{P}$ and $\mathcal{N}$ have their constant terms in any pair of
complementary subalgebras of $\gtilde_0$, the subalgebra of $\gtilde$
which consists of constant
matrices. See  Section \ref{scheme},
and in particular Remark \ref{parallelremark}, 
for a specific instance of this.

We will construct flows on the finite dimensional vector subspaces
\bdm
\gtilde_{d}^1 := \{ X = \sum_{i=d}^1 X_i z^i \in \gtilde \},
\edm
where $d$ is any integer less than or equal to $1$.
Therefore we can assume that all
maps which are to be differentiated are between finite-dimensional
subspaces of the infinite-dimensional loop algebras.

 Let 
\bdm
\gtilde^1 := \bigcup_{d \leq 1} \gtilde_d^1.
\edm

\begin{proposition}  \label{loopgroupproposition}
Let $\pip: \gtilde \to \mathcal{P}$ be the projection.
Let $p_1$,..., $p_n \in {\mathbb C}[x,y]$ be a collection of polynomials 
in two variables. Let $V_i: \, \gltilde \to \gltilde$ be the map defined
by 
\bdm
V_i(X) := p_i( X , z^{-1}),
\edm
and suppose that, for all $i$,
\beq \label{conditionA}
V_i(\gtilde^1) \subset  \gtilde^1,
\eeq
i.e.  $p_i (X,z^{-1})$ is in $\gtilde^1$ whenever $X$ is. 

 Then: 
\begin{enumerate}
\item The Lax equation
\beq  \label{flowequation}
\dd X = [X, \sum_{i=1}^n \pip V_i(X) \dd t_i],
\eeq
on $\gtilde^1$, is of \emph{finite type}: if 
$X \in \gtilde_d^1$ then so are the coefficients of the 1-form $\dd X$.
\item
Equation (\ref{flowequation}) is Frobenius integrable, i.e. $\dd ^2 X =0$.
\item
The 1-form $A:= \sum_{i=1}^n \pip V_i(X) \dd t_i$ satisfies the 
\emph{Maurer-Cartan equation}
\bdm
\dd A + A \wedge A = 0.
\edm
\end{enumerate}
\end{proposition}

\begin{definition}
We will call the equation (\ref{flowequation}) for the case
$\gtilde = \gtildesig$ a \emph{$k$-symmetric AKS system}.
\end{definition}

\begin{remark}  \label{conditionremark} 
 The condition (\ref{conditionA})  is merely to ensure that the equation
(\ref{flowequation}) makes sense on $\gtilde^1$. Finding polynomials
which map $\gtilde \to \gtilde$ depends on the Lie algebra, but ensuring
that the highest power of $z$ in $p_i(X,z^{-1})$ is 1 is easily achieved
by multiplying by a sufficiently high power of $z^{-1}$.  

  We note here
the useful fact that for the above-mentioned example $\gtilde _{\sigma}$,
when $\sigma$ is an order $k$ inner automorphism $\sigma (X) = QXQ^{-1}$ for
some $Q \in G$ (a Lie group with Lie algebra  $\mathcal{G}$),
if the monomial $p(X) = X^j$ takes $\mathcal{G} \to \mathcal{G}$ then
$p$ also takes $\gtilde _{\sigma} \to \gtilde _{\sigma}$, and 
so does $z^{lk}p$ for any integer $l$.  
\end{remark}

 \begin{remark}  
Given a Lax 
equation $\dd X = [X,A]$, our choice of $A$ is not unique.
For example we can always add to $A$ any 1-form whose coefficients commute
with $X$, without changing the equation.
The fact that, for $k$-symmetric AKS systems, $A$ satisfies
 the Maurer-Cartan equation says something about our choice: 
by standard Lie group theory,  the satisfaction of this equation
is equivalent (on a simply connected domain) to the fact
that there exists $F:  {\mathbb R}^n \to \widetilde{\textup{GL}}(n,{\mathbb C})$
 such that $A = F^{-1} \dd F$, $F(0) = I$ (here $F^{-1}$ denotes the
matrix inverse). This give us the useful formula
\bdm
X(t) = F^{-1}(t) X(0)F(t),
\edm
which is verified by checking that the right hand side solves the
Lax equation $\dd X = [X, \, F^{-1} \, \dd F]$, which defines $X$.
 \end{remark}

\subsection{Proof of Proposition \ref{loopgroupproposition}} 
The fact that the flow is well-defined on $\gtilde_d^1$ is easily
verified: if $X \in \gtilde_d^1$ then  the condition 
(\ref{conditionA}) ensures that $\pip V_i (X)$ are in $\gtilde$, and so,
 therefore, are the coefficients of $\dd X$, since $\gtilde$ is a 
subalgebra.  Now if 
$X = ( ... + X_0 + X_1 z)$ then (\ref{conditionA}) also implies 
that $V_i(X) = ... + C_1 X_1^k z$, for some non-negative integer
$k$. Thus $\pip V_i(X) = C_0 + C_1 X_1^k z$, where $C_0$ and
$C_1$ are constants, and 
\beqas
[X, \, \pip V_i (X)] &=& [ \sum_{-d}^1 X_i z^i , \, C_0 + C_1 X_1^k z ] \\
&=& [ \sum_{-d}^1 X_i z^i , \, C_0] +
  [ \sum_{-d}^0 X_i z^i , \, C_1 X_1^k z ] + [X_1 z, C_1 X_1^k z ].
\eeqas
The last term is zero, so it follows that $\dd X$ takes it's values in $\gtilde_d^1$.

To prove parts 2 and 3, we will use some standard theory for
 constructing flows from
ad-invariant functions on Lie groups, which we now recall.

 Let $\mathcal{G}$ be any Lie algebra.
\begin{definition} An \emph{ad-equivariant vector field} is a map $V:~\mathcal{G}
\to \mathcal{G}$ which satisfies
\beq \label{adinvariant}
\dd V \Big |_X ([X,Y])  = [V(X), Y].  
\eeq
\end{definition}
\begin{remark} 
If $\mathcal{G}$ is equipped with an invariant inner product, then 
ad-equivariant vector fields are  (locally) gradients,
with respect to this inner product, of ad-invariant polynomials on
$\mathcal{G}$.  Thus, for example, if $\mathcal{G}$ is the Lie algebra of
a compact semi-simple
Lie group of rank $k$ then we can expect to find at most $k$ 
polynomials which are linearly independent pointwise.
 This follows from the fact that an ad-invariant function is given by
its values on a maximal torus, which has dimension $k$.
\end{remark}

 For the purpose of proving Proposition \ref{loopgroupproposition} 
we note the following:
\begin{lemma}
Any polynomial $p(X,z,z^{-1})$  in $X$, $z$ and $z^{-1}$ is an ad-equivariant
vector field on $\widetilde{\textup{gl}}(m,{\mathbb C})$. 
\end{lemma}
 \begin{proof}  This is easily verified by first checking that
$V(X) = X^n$ satisfies the equation (\ref{adinvariant}) as follows:
\bdm
V(X + t[X,Y]) = X^k + t(\sum_{i=0}^{n-1}X^{n-1-i}(XY-YX)X^i) + o(t^2),
\edm 
so that
\beqas 
\dd V \Big |_X ([X,Y]) &=& X^n Y - Y X^n \\
 &=& [X^n, Y].
\eeqas
Hence the observation that both the left and right hand sides of (\ref{adinvariant})
are linear over ${\mathbb C} [z,z^{-1}]$ proves the lemma. 
\end{proof}

 Now the proof of Proposition \ref{loopgroupproposition} follows 
from the following standard result 
(see, for example, \cite{burstallpedit1994}):

\begin{proposition}  \label{integrableflows}
Suppose we have a collection of $n$ ad-equivariant vector fields $V_i$ on 
a Lie algebra $\mathcal{G} = \mathcal{P} \oplus \mathcal{N}$, where
$\mathcal{P}$ and $\mathcal{N}$ are both subalgebras,
 and consider the projections
\bdm  
A_i := \pi _{\mathcal{P}} V_i.
\edm
Then: 
\begin{enumerate}
\item The system 
\beq  \label{flowequations}
X_{t_i} = [X, A_i(X)], ~~~~  i= 1...n,
\eeq
on $\mathcal{G}$ is Frobenius integrable.
\item If $X: {\mathbb R}^n \supset U \to \mathcal{G}$ is a solution to 
(\ref{flowequations}) then the 1-form
\bdm
A := \sum_{i=1}^{n} A_i \dd t_i
\edm
satisfies the Maurer-Cartan equation $dA + A \wedge A = 0$.
\end{enumerate}
\end{proposition}

\begin{proof} 
Using the Jacobi identity, it is straightforward to reduce the 
first statement to the equation:
\bdm
[X,  [A_i,A_j] + \nabla_{[X, A_i]} A_j - 
\nabla_{[X, A_j]} A_i ] = 0 ~~~~~~\textup{for all} ~i,j.
\edm 
Additionally, the  Maurer-Cartan equation of statement two is equivalent to
\bdm
[A_i, A_j] +
\frac{\partial}{\partial t_i} A_j - \frac{\partial}{\partial t_j} A_i 
 = 0   ~~~~~~ \textup{for all}~i,j.
\edm
If the first statement holds, then the Maurer-Cartan equations can be 
rewritten as
\beq  \label{mcequations}
[A_i,A_j] + \nabla_{[X, A_i]} A_j - 
\nabla_{[X, A_j]} A_i = 0 ~~~~~~\textup{for all} ~i,j.
\eeq
Thus it suffices to verify the equations (\ref{mcequations}). This
follows from Lemma \ref{mclemma} below.
\end{proof}

We state here a result which is more general than needed, but is of
interest as it can be used to define variants of the AKS theory,
as in \cite{feruspedit1996}.
\begin{lemma} \label{mclemma}
Suppose $\mathcal{G} = \mathcal{P} \oplus \mathcal{N}$ is a decomposition of a 
Lie algebra into vector subspaces, and that $V_i$ are ad-equivariant vector fields
on $\mathcal{G}$. Let $A_i := \pi _{\mathcal{P}} V_i$.  Then:
\begin{enumerate}
\item
\bdm
 [A_i,A_j] + \nabla_{[X, A_i]} A_j - 
\nabla_{[X, A_j]} A_i = 
[\pip V_i, \pip V_j]
+ \pip \{ [V_j, \pip V_i ] - [V_i, \pip V_j] \}. \\
\edm
\item If $\mathcal{P}$ is  Lie sub-algebra then
\bdm
 [A_i,A_j] + \nabla_{[X, A_i]} A_j - 
\nabla_{[X, A_j]} A_i 
= \pip [\pin V_i, \pin V_j].
\edm
\end{enumerate}
\end{lemma}
 \begin{proof}
In both of the above equations, the left hand side is equal to
\bdm
 [\pip V_i, \pip V_j] + 
\nabla_{[X, \pip V_i]} \pip V_j - \nabla _{[X, \pip V_j]} \pip V_i. 
\edm
Using the decomposition of $\mathcal{G}$, this amounts to
\beqas
[\pip V_i, \pip V_j] + \pip \{ \nabla_{[X, \pip V_i]} V_j - 
  \nabla _{[X, \pip V_i ]} \pin V_j - \nabla_{[X, \pip V_j]} V_i 
  + \nabla_{[X, \pip V_j]} \pin V_i \} \\
~~ = [\pip V_i, \pip V_j] + \pip \{ \nabla_{[X, \pip V_i]} V_j
 - \nabla _{[X, \pip V_j]} V_i \}.
\eeqas
Using the ad-equivariance of $V_i$ and $V_j$, this equates to the right hand
side of the equation in statement 1.

To get statement 2, we use the assumption that $\mathcal{P}$ is a subalgebra
to obtain:
\beqas
  \pip \{ [V_j, \pip V_i ] - [V_i, \pip V_j] \} \\
~~ = 
[\pip V_j, \pip V_i ] + \pip [\pin V_j, \pip V_i ] + \pip [\pip V_j, \pip V_i ]
+ \pip [\pip V_j , \pin V_i] \\
~~ = -[\pip V_i, \pip V_j] + \pip  [V_j, V_i] - \pip [\pin V_j, \pin V_i] \\
~~ = -[\pip V_i, \pip V_j] + \pip [\pin V_i, \pin V_j].
\eeqas
Here we used the easily verified fact that ad-equivariant vector fields commute.
Now substitute this into statement 1 to get statement 2.
\end{proof}

\subsection{Linearization of a System of Finite Type}  \label{linearization}

 Recall that to describe a 
holomorphic line bundle over ${\mathbb CP}^1$ one merely needs to 
provide the transition function between the local trivializations
over $\mathbb C$ and ${\mathbb C}^* \cup \{ \infty \} $ respectively.
Similarly,  since $ z : \Sigma \to {\mathbb CP}^1 $ is a smooth covering of 
${\mathbb CP}^1$, it follows that $\Sigma ^+ = 
\Sigma - z^{-1}(\infty)$ and $\Sigma^- = \Sigma - z^{-1}(0)$
is an open cover of $\Sigma$ by two sets which admit only trivial
holomorphic line bundles. Thus 
a holomorphic line bundle over $\Sigma$ is defined by a transition
function between local trivializations over $\Sigma ^+$ and 
 $\Sigma^-$, that is,  a  holomorphic function $\rho : 
\Sigma^+ \cap \Sigma ^- \to {\mathbb C}^*$. More specifically, if 
$s^{\pm} (z, w, t)$ is a non-vanishing section of  $L(t)$ over
$\Sigma ^{\pm}$, then $\rho$ must satisfy the equation $s^+ = \rho s^-$ on
$\Sigma^+ \cap \Sigma^-$.\\
  
 The following proposition,
 says that the evolution of $L(t)$ in the moduli space
of line bundles over $\Sigma$ is linear.

\begin{proposition} \label{linearizationprop}
Suppose that $X(z, t)$ is a solution of the system (\ref{flowequation})
with the initial condition $X(z, 0) = X_0(z)$. Let  $\Sigma$
be the spectral curve of $X_0$ and $\rho_0$ the transition function of
the eigenspace bundle $L_0$.  Let $\mu_i : \Sigma \to {\mathbb C}$ be
the eigenvalue
function of the operator $V_i(X_0(z))$ defined by the equation
\bdm
V_i((X_0(z))) s_i(z, w) = \mu_i(z, h) s_i(z, w),
\edm
where $w$ is the eigenvalue function of $X_0(z)$. Denote
by $\bf \mu$ the $n$-tuple $[\mu _1,...,\mu _n]$.

 Then the eigenspace bundle $L(t)$ of 
$X(z, t)$ is defined by the transition function   
\bdm
\rho (t) = \rho _0 \exp(- {\bf t \cdot \mu}).
\edm
\end{proposition}
 \begin{proof} 
The proof given here is based on an argument given in \cite{fpps}.
  Let $F(z, t)$ be a solution of 
\bdm
F^{-1}\dd F = \sum_{i=1}^n \pip V_i \dd t_i, ~~~~ F(z, 0) = I, ~~z \in {\mathbb C}.
\edm
Since $\pip V_i$ are polynomials in $z$, $F$ is a holomorphic function of 
$z$ on ${\mathbb C}$.
  As described earlier, we know that $X(t) = F^{-1}(t) X_0 F(t)$,
so we immediately have $X(t) F^{-1}(t)  s^+(0) = w F^{-1}(t) s^+(0)$
 (given that 
$s^+(0)$ is an eigenvector of $X_0)$.  In other words, if $s^+(0)$ is non-vanishing
section of $L_0$ over $\Sigma^+$, then 
\bdm
s^+(t) := F^{-1}(t) s^+(0)
\edm
is a non-vanishing section of $L(t)$ over $\Sigma^+$. 
 
  We need a similar formula for $s^-$ over $\Sigma^-$.
Set
\bdm
G := \exp(- \sum_{i=1}^n t_i V_i(X_0) F(t)).
\edm
Then
\beqas
G^{-1} \dd G &=& F^{-1}[\sum_{i=1}^n -V_i(X_0) F \dd t_i + \dd F] \\
&=&   \sum_{i=1}^n[ - F^{-1} V_i(X_0) F + \pip V_i (X(t))] \dd t_i \\
&=& \sum_{i=1}^n [-V_i(X(t)) + \pip V_i(X(t))] \dd t_i \\  
&=& \sum _{i=1}^n - \pin V_i(X(t)) \dd t_i.
\eeqas
We used the fact that $V_i(X)$ is of the form $p_i(X,z^{-1})$ and that
$X(t) = F^{-1} X_0 F$ on the third line.

 Since $\pin V_i (X)$ is 
polynomial in $z^{-1}$, it follows that $G$ is holomorphic on
${\mathbb C}^* \cup \{ \infty \}$.  Using the fact that $V_i(X_0)$ commutes
with $X_0$, we also  obtain the formula $X(t) = G^{-1}(t)X_0G(t)$, and
therefore, as with $s^+$,
\bdm
s^-(t) = G^{-1}(t) s^-(0)
\edm
is a non-vanishing section of $L(t)$ over $\Sigma^-$. 

 Hence, on the intersection $\Sigma^+ \cup \Sigma^-$,
 one obtains
 \beqas 
s^+(t) & =& F^{-1} \rho_0 s^-(0) \\
&=& \rho_0 G^{-1} G F^{-1} s^- (0) \\
&=& \rho_0 G^{-1} \exp(-\sum_{i=1}^n t_i V_i(X_0)) s^-(0) \\
&=& \rho_0 G^{-1} \exp(- \sum_{i=1}^n t_i \mu_i) s^-(0) \\
&=& \rho_0 \exp(-{\bf t \cdot \mu}) s^- (t). 
\eeqas 
\end{proof}

\begin{remark} \textup{ It is essential to the argument that
the decomposition of the loop algebra $\mathcal{G} = \mathcal{P}
\oplus \mathcal{N}$ is into subspaces consisting of functions which
are holomorphic in $z$
and $z^{-1}$ respectively, and that the functions $V_i(X)$ are
polynomials in $X$ whose coefficients are meromorphic functions of
$z$. Any such commuting flows on loop
algebras will linearize as described here. } \end{remark}


\section{Simple $k$-symmetric systems and their properties.} 
 \label{properties}

In this section we are interested in the question of periodicity of
solutions. 
 Since the evolution of the associated line bundle 
is essentially linear for $k$-symmetric AKS systems, its periodicity
can be characterized in terms of the initial data 
 (see \cite{fpps}). Clearly periodicity of $F$ implies 
periodicity of the line bundle. Below we will
look at what periodicity of the line bundle implies about the 
corresponding solution $F$.

Throughout this section we will be studying the solutions 
of a restricted class of AKS systems.

\begin{definition} A $k$-symmetric AKS system will be called
\emph{simple} if the decomposition (\ref{decomposition}) is
of the form
\bdm
\gtildesig = \pa \oplus \na,
\edm
where $\pa$ and $\na$ are defined as in (\ref{simpledecomp}),
so that  the whole subalgebra of constants, $\mathcal{V}_0$,
 is contained in $\na$. 
\end{definition}

For simple systems, the projection of $\gtildesig^1$ onto $\pa$
is just $\pipa (... + X_0 + X_1 z) = X_1 z$.  This allows us to
prove the results below more easily.

\subsection{Recovering $F(t)$ from the Spectral Data} 

Recall that if $F(0) = I$ then we have the formula
\beq  \label{magicformula}
X(z,t) = F^{-1}(t) X(0) F(t).
\eeq
For a more general initial condition $F(0) = F_0 \in \gtilde$, set 
\bdm
Y_0 := F_0 X(0) F_0^{-1},
\edm
 and the following formula is similarly deduced:
\bdm
 X(t) = F^{-1}(t) Y_0 F(t).
\edm

The following proposition
says that for simple $k$-symmetric AKS systems, the spectral data 
determine $X$ and $F$ uniquely from their initial conditions.
\begin{proposition} \label{uniquenessprop}
Suppose $X(z,t) = F^{-1}(z,t) Y_0(z) F(z,t)$ is a solution
corresponding to the regular initial conditions
 $X(z,0) = F_0^{-1} (z) Y_0(z) F_0(z)$ and $F(z,0) = F_0(z)$,
obtained from Proposition \ref{solutionproposition} 
for the loop algebra decomposition $\mathcal{G} = \mathcal{P}_a \oplus \mathcal{N}_a$.
Then both $X(z,t)$ and $F(z,t)$ are determined by the spectral data 
$\{ \Sigma, L(t), z \}$, the initial conditions for $X$ and $F$, and the
requirement that $F^{-1} \dd F \in (\gtildesig)_1^1 \otimes \Omega ({\mathbb R}^n)$.  
\end{proposition}

\begin{proof} 
 As described in Section \ref{reconstruct},
$X(z,t)$ is determined by its spectral data up to conjugation by an 
element of  $\textup{GL}(2n,{\mathbb C})$, so suppose we have
\beqas
\hat{X}(z,t) &=& B^{-1}(t) X(z,t) B(t) \\
&=& B^{-1}(t) F^{-1}(z,t) Y_0(z) F(z,t) B(t),
\eeqas
which has the spectral data $\{ \Sigma, L(t), z \}$ and the given initial
conditions, that is,
$\hat{X} (z,0) = F_0^{-1}(z) Y_0(z) F_0(z)$ and 
suppose that $B(t)$ is chosen so that
$F(z,0)B(0) = F_0(z)$, and such that
$(FB)^{-1} \dd (FB) \in (\gtildesig)_1^1 \otimes \Omega ({\mathbb R}^n)$.
Note that these assumptions mean that $B(0) = \textup{Id}$.
  We show that 
$B(t) = \textup{Id} = \textup{constant}$, which proves that $\hat{X} = X$.
Since $F$ is determined by $X(t)$ and $F(0)$ through the equation
$F^{-1} \dd F =  \sum_{i=1}^n \pip \mathcal{V}_i(X) \dd t_i$, this proves
the proposition.

 The requirement that  $F^{-1} \dd F \in (\gtildesig)_1^1 \otimes 
\Omega ({\mathbb R}^n)$ simply means we can write
 $F^{-1} \dd F = A_1z$, where $A_1 \in V_1 \otimes \Omega({\mathbb R}^n)$, and
then 
\beqas
(FB)^{-1} \dd (FB) &=& B^{-1}F^{-1} [(\dd F) B + F \dd B] \\
&=& B^{-1} A_1 B z + B^{-1} \dd B.
\eeqas
$B$ has no $z$ dependence, so for the right hand side to lie in
$(\gtildesig)_1^1 \otimes \Omega ({\mathbb R}^n)$
 we must have $B^{-1} \dd B = 0$, which implies that
$\dd B =0$ and $B= \textup{constant} = \textup{Id}$.
\end{proof}

\subsection{Periodicity of Solutions}
In this section we study the implications of periodicity of the
spectral data in the case of simple $k$-symmetric systems.


\subsubsection{A Translation Lemma}
Every initial condition $X_0$ in
$(\gtildesig)^1$ generates a polynomial Killing field $X(t)$
 corresponding to a solution $F(t)$ with initial condition $F(0) = I$.
It is of interest to know how $F$ is related to the solution 
corresponding to
that generated by the initial condition $\hat{X}(0) = X(Q)$, for some 
$Q \in {\mathbb R}^n$. 
\begin{lemma} \label{translationlemma}
Let $F$ be a simple $k$-symmetric map with Killing field $X$,
with the initial conditions $X(0) = X_0$ and $F(0) = I$.
Let $\hat{X}$ and $\hat{F}$ be the solutions corresponding to the initial
conditions $\hat{X}(0) = X(Q)$ and $\hat{F}(0) = I$. Then, for all 
$t \in {\mathbb R}^n$, 
\beq
F(t+ Q) = F(Q) \hat{F}(t).
\eeq
\end{lemma}

\begin{proof} 
Since $Y(t):= X(t+Q)$ is a solution of
the Lax equation $\dd Y = [Y, \, A(Y)]$, and $Y(0) = X(Q)$, it follows
by uniqueness that $Y = \hat{X}$, which is to say that 
\beqas  \label{xhat}
\hat{X}(t) &=& X(t + Q).
\eeqas
We want to show that 
\bdm
\hat{F}(t) = F^{-1}(Q) F(t+Q).
\edm
It is enough to check that the right hand side satisfies the equations
defining $\hat{F}$, namely $\hat{F}(0) = I$ and
 $\hat{F}^{-1} \dd \hat{F} = A(\hat{X}) :=  
  \sum_{i=1}^n \pip \mathcal{V}_i(\hat{X}) \dd t_i$. The first equation 
clearly holds, and the second can be verified as follows:
\beqas
F^{-1}(t+Q)F(Q)F^{-1}(Q) \dd F(t+Q) &=& F^{-1}(t+Q) \dd F(t+Q) \\
&=& A(X(t+Q)) \\
&=& A(\hat{X}).
\eeqas
\end{proof}

 A consequence of Lemma \ref{translationlemma} is that 
if $F$ is not injective then it has a period:

\begin{lemma} \label{lemma4}
Let $F$ be a simple $k$-symmetric map. If $F(a) = F(b)$ for some
$a \neq b$ then $b-a$ is a period for $F$.
\end{lemma}
\noindent \textbf{Proof:} $~~$ We first assume that $a=0$ and $F(0)= F(b) = I$.
If $X$ is a Killing field for $F$, then $X(b) = F^{-1}(b)X(0)F(b) = X(0)$,
so it follows immediately from Lemma \ref{translationlemma}
 that $F(t+b) = F(t)$.
The general case can be deduced from this one by setting 
$\hat{F}(t) := F(a)^{-1} F(t+a)$, so that $\hat{F}(0) = \hat{F}(b-a) = I$.

\subsubsection{Periodicity of the Killing Field}

\begin{definition}
A map $F$ from ${\mathbb R}^n$ to a Lie group $G$ is \emph{type I quasiperiodic}
with quasiperiod $P \in {\mathbb R}^n$ if 
\bdm
F(t+P) = F_0 \, F(t)
\edm
 for all $t \in {\mathbb R}^n$, where $F_0$ is some
constant element of $G$.
\end{definition}
The translation Lemma of above can be used to 
show that periodicity for $X$ is equivalent to type I quasiperiodicity for $F$.
In fact we only need $X(P) = X(0)$ to guarantee both of these properties:

\begin{proposition} \label{killingprop}
Let $F$ be a simple $k$-symmetric map with initial condition
 $F(0) = I$ and Killing field $X$. 
The following three  conditions are equivalent:
\begin{enumerate}
\item 
  $X(P) = X(0)$.
\item
  $F$ is type I quasiperiodic with quasiperiod $P$, and $X(P) = X(0)$.
\item
  $X(t+P) = X(t)$ for all $t \in {\mathbb R}^n$.
\end{enumerate}
\end{proposition}
\begin{proof} 
We prove $(1) \Rightarrow (2) \Rightarrow (3) 
\Rightarrow (1)$. \\

\noindent ${\bf (1) \Rightarrow (2):}$ \\
This follows from 
Lemma \ref{translationlemma}. 

\noindent ${\bf (2) \Rightarrow (3):}$ \\ 
If $F(t+P) = F_0 F(t)$ then 
$F^{-1} \dd F \Big |_{t+P} = F^{-1} \dd F \Big |_t$. But recall that
$F^{-1} \dd F = \sum_{i=1}^n \pipa V_i(X) \dd t_i,$
 where, if
\bdm 
X= X_{-d} z^{-d} + ... + X_0 + X_1 z,
\edm
 then $\pipa V_i(X) := X_1^{2i-1} z.$ In other words, we have the formula
\bdm
F^{-1} \dd F = \sum_{i=1}^n X_1^{2i-1} z \dd t_i.
\edm
Thus if the connection $F^{-1} \dd F$ is periodic, then so, certainly, is the function
$X_1$.  If we now consider the differential equations satisfied by  $X$,
namely (\ref{flowequation}), 
we see this reduces to
\beqa
\frac{\partial X_i}{\partial t_j} &=& [X_{i-1}, \, X_1^{2j-1}],
 ~~~~ i = (-d+1),...,1,   \label{de1} \\
\frac{\partial X_{-d}}{\partial t_j} &=& 0.
\eeqa
Thus $X_{-d}$ is constant, and so $X_1$ and $X_{-d}$ both have a 
period $P$. Hence by the equation (\ref{de1}), any partial derivative of
$X_{-d+1}$ also has a period $P$, and so, therefore, does  $X_{-d +1}$.
Similarly, by induction,
$X_i$ has a period $P$ for all $i$. \\

\noindent ${\bf (3) \Rightarrow (1):} $ \\
Immediate.
\end{proof}


\subsubsection{Periodicity of the Eigenspace Bundle}
We can now proceed to the main goal of this section,  the proof of Proposition \ref{proposition1}.
\begin{definition}
Let $G$ be a subgroup of $\gltilde$. 
A map $F$ from ${\mathbb R}^n$ to $G$ is \emph{type II quasiperiodic}
with quasiperiod $P \in {\mathbb R}^n$ if 
\bdm
F(t+P) = F_0 \, B^{-1} F(t) B,
\edm
 for all $t \in {\mathbb R}^n$, where $F_0$ is some
constant element of $G$, and $B \in \textup{GL}(m, {\mathbb C})$.
\end{definition}

\begin{proposition} \label{periodprop}
Suppose that $F$ is a simple $k$-symmetric map with Killing field
$X$ and  the initial conditions $X(0) = X_0$ and $F(0) = I$.
Suppose further that the eigenspace bundle $L(t)$ of $X(t)$ satisfies
 $L(t+P) = L(t)$ for all $t$.  Then 
 $F$ is type II quasiperiodic with quasiperiod $P$, in particular, 
there exists a matrix $B \in \textup{GL}(2n, {\mathbb C})$ such that
we have for all $t$:
\beqas
F(t+P) = F(P) B^{-1}F(t)B.
\eeqas
\end{proposition}
Before proving this, 
 we will find out what can be said concerning
the Killing field $X$ under these circumstances.
\begin{lemma} \label{lemma1} 
If the hypotheses of Proposition \ref{periodprop} hold then
\beq
X(t+P) = B^{-1}X(t)B,
\eeq
for some constant matrix $B \in \textup{GL}(2n, {\mathbb C})$.
\end{lemma}
\begin{proof}  
 As described in Section \ref{reconstruct},
$X(z,t)$ is determined by its spectral data up to conjugation by an 
element of  $\textup{GL}(2n, {\mathbb C})$. So we know that there exists a
matrix $B(t)$ such that 
\bdm
X(t+P) = B^{-1}(t) X(t) B(t).
\edm
We know that $Y(t) := X(t+P)$ is the unique solution of the differential equation:
\beqas
\dd Y = [Y, \, A(Y)],\\
Y(0) = X(P).
\eeqas
We show that $\hat{Y}(t) := B^{-1}(0) \, X(t) \, B(0)$ also satisfies these equations,
which will complete the proof. Clearly the initial condition $\hat{Y}(0) = X(P)$
holds, so let us check the Lax equation:
\beqas
\dd \hat{Y} &=& B^{-1}(0) \dd X B(0)  \\
  &=& B^{-1}(0) \, [X, A(X)] \, B(0) \\
  &=& [B^{-1}(0) \, X \, B(0), \, B^{-1}(0) \, A(X) \, B(0)] \\
  &=& [\hat{Y}, \, B^{-1}(0) \sum_{i} \pipa V_i(X) \, B(0)]\\
  &=& [\hat{Y}, \, \sum_i \pipa V_i (B^{-1}(0) \, X \, B(0)) ] \\
  &=& [\hat{Y}, \, A(\hat{Y})].
\eeqas
The step on the fourth line was allowed because $B$ does not depend on the
parameter $z$, and therefore conjugation by $B$ commutes with projection
onto $\mathcal{P}_a$. 
\end{proof}

\begin{remark} 
The above proof depends on the fact that, for simple $k$-symmetric
systems, our decomposition
$\gtilde = \pa \oplus \na$ is such that all the constant
terms are in just one of $\pa$ or $\na$. Otherwise 
conjugation by a constant matrix $B$ will in general
 not commute with projection
onto $\mathcal{P}$. \end{remark}

\begin{proof}{(of Proposition \ref{periodprop}).} 
The equation
\beqas
F(t+P) = F(P) B^{-1}F(t)B
\eeqas
can be proved by checking  that both the left
and right hand sides  are solutions,
$H$, for the differential equation:
\beqa
H^{-1} \dd H = B^{-1} F^{-1} \,  \dd F \, B,  \label{difeq}\\
H(0) = F(P). \label{intcon}
\eeqa
The initial condition (\ref{intcon}) evidently holds, so we verify
 the differential 
equation (\ref{difeq}) here: setting
\bdm
 H(t) = F(t+P)
\edm
we have
\beqas
H^{-1} \dd H &=& F^{-1}(t+P) \dd F (t+P)  \\
 &=& A(X(t+P)) \\
 &=& A(B^{-1} X(t) B) \\
 &=& B^{-1}A(X(t))B \\
 &=& B^{-1}F^{-1} \dd F B. 
\eeqas
As in the proof of Lemma \ref{lemma1} we again used the fact that 
conjugation by $B$ commutes with projection onto $\mathcal{P}_a$ on
the fourth line. 

On the other hand, setting
\bdm
H(t) = F(P)B^{-1}F(t)B,
\edm
one obtains:
\beqas
H^{-1} \dd H &=& B^{-1}F^{-1}(t)BF^{-1}(P)F(P)B^{-1} \dd F(t) B \\
 &=& B^{-1} F^{-1} \dd F B.
\eeqas
\end{proof}


\section{Finite Type Flat Immersions in a Sphere} \label{applications}
 
\subsection{Admissible Connections} \label{sectadmissconn}
Suppose $f:  {\mathbb{R}^n} \supset U\to S^{2n-1} \subset {\mathbb{R}}^{2n}$ is an 
immersion. One can associate to $f$ an  immersion $F: U \to SO(2n)$,
(called an \emph{adapted framing}) of the form:
\bdm
 F = [e_1,...,e_n,\xi _1,...,\xi _n],
\edm
where $\xi _1 = f$ and the vectors $\{e_i\}$ and
 $\{\xi _2, ..., \xi _n \}$ are (non-unique) orthonormal bases for the tangent and normal
spaces respectively to $f(U) \subset S^{2n-1}$. 

 Denote by $A= F^{-1}\dd F
 \in so(2n) \otimes \Omega^1(U)$ the pull-back by $F$ of the left 
Maurer-Cartan form of $SO(2n)$. More explicitly,
\bdm  
A  = \bbar {cc}
      \omega & - \beta ^T \\
     \beta  & \eta \\
  \ebar, \label{adaptedform}
\edm
where $\omega \in so(n) \otimes \Omega^1(U)$ is the
Levi-Civita connection form on $f(U)$ and $\eta \in so(n) \otimes 
\Omega^1(U)$
and $\beta$ are the normal connection and the second fundamental form 
respectively of the immersion $f$ into ${\mathbb{R}}^{2n}$. 

 Since $\xi _1 = f$,
the requirement that $\{ e_i \}$ span the tangent space to $f$ is equivalent to
the statement that the first row and first column of $\eta$ consist of
zeros.  The map $f$ will be immersive provided that the first row of
the second fundamental form $\beta$ (which is the dual frame to $\{e_i \}$)
consists of $n$ linearly independent 1-forms.

 Given that $A = F^{-1} \dd F$, we can calculate 
the \emph{Maurer-Cartan equation} 
\bdm 
\dd A + A \wedge A = 0.
\edm
For 1-forms $A$ which arise in the way described above, this is equivalent 
to the structure equations and the Gauss, Codazzi
and Ricci equations for the immersion $f$. 

Conversely, by standard theory of Lie groups, if $A \in so(2n) \otimes \Omega^1(U)$ 
satisfies the Maurer-Cartan
equation then, provided $U$ is simply connected,
 $A$ can be integrated to give a map $F: U \to SO(2n)$. 

 Now suppose we have a family $A^{z}$ of such 1-forms, of the form:
\beq
A^{z} = \bbar {cc}
                \omega &  - z \beta ^T \\
                z \beta  & \eta \\
	      \ebar,
\eeq
such that the first row and column of $\eta$ are zero and such
that the Maurer-Cartan equation $\dd A^{z} + A^{z} 
\wedge A^{z} = 0$ is satisfied for all $z \in  {\mathbb{C}}$, and
such that $A^z$ is real (i.e. the components of the matrix are
real valued one-forms) for $z \in {\mathbb R}$.
 Then for each non-zero real value of
$z$  we can locally integrate $A^{z}$ to obtain
an adapted frame $F^{z}$ for a map $f^{z}:
M \to S^{2n-1}$. Moreover, writing out the Maurer-Cartan equations
explicitly,
\beqas
0 = 
\bbar {cc}
      \dd \omega & -z \dd \beta^T \\
     z \dd \beta & d \eta \\
\ebar +
\bbar {cc}
    \omega \wedge \omega - z^2 \beta^T \wedge \beta &
        z (\omega \wedge \beta - \beta ^T \wedge \eta) \\
    z(\beta \wedge \omega + \eta \wedge \beta) & 
        -z^2 \beta \wedge \beta ^T + \eta \wedge \eta 
\ebar,
\eeqas
and equating coefficients of like powers of $z$ we obtain an
additional condition:
\bdm 
\dd \omega + \omega \wedge \omega = 0,
\edm
which says that the curvature of the induced metric on $f(U) \subset S^{2n-1}$
is zero. 

If $f^z$ are immersions (this doesn't depend on $z$), then
the equations obtained from the other components give nothing
further, since all of them follow from the Maurer-Cartan equation at
a single value of $z$, with the exception of the flatness of 
the normal bundle, $\dd \eta + \eta \wedge \eta = 0$; but it is well known
that for any (local) isometric immersion of flat space $E^n$ into $S^{2n-1}$ 
the normal bundle is necessarily flat. 

 Now given the flatness of the normal bundle, we could choose
a parallel frame $\{ \xi_1, ...,\xi_n \}$ for the normal bundle,
 i.e. one in which $\eta =0$.  This leads to the following:
\begin{definition}
An \emph{admissible connection} on ${ \mathbb R}^n$ is a family of 
so(2n)-valued 1-forms $A^{z}$ of the form
\beq  \label{admissible}
A^{z} = \bbar {cc}
                \omega &  - z \beta ^T \\
                z \beta  & 0 \\
	      \ebar,
\eeq
(where $\omega$ is $n \times n$) which satisfy the Maurer-Cartan equation
for all $z \in {\mathbb C}$.
\end{definition}
Note that an admissible connection can be integrated to get a family
of maps $F^z: {\mathbb R}^n \to SO(2n)$, $z \in {\mathbb R} _*$,
and any of the columns $(n+1),...,2n$ of $F$ is a flat immersion into
$S^{2n-1}$, provided that it is injective.    Conversely,
any flat immersion of ${\mathbb R}^n$ into $S^{2n-1}$ gives rise to an admissible
connection. Therefore, we now consider the problem of producing
admissible connections.

\begin{remark} \textup{We should note here that in other applications of 
AKS theory to special submanifolds, it is always the case
that one must be able to insert a spectral parameter $z$ into the 
relevant Maurer-Cartan equation in some such way as we have done here.
However, it is not necessary that the satisfaction of the equation for
all $z$ need add anything new, as it did here by introducing flatness. 
 See, for example, \cite{fpps}.} \end{remark}

\subsection{A Scheme for Generating Admissible Connections}   \label{scheme}

Let $so(2n,{\mathbb C})$ and $\sigma = Ad(Q)$ be as in the introduction,
so that $so(2n,{\mathbb C}) = \mathcal{V}_0 \oplus \mathcal{V}_1$,
where $\mathcal{V}_i$ is the eigenspace of $\sigma$ corresponding to $(-1)^i$. 
More explicitly,
\bdm
\mathcal{V}_0 = \left\{ \bbar  {cc}
	       {*} & 0 \\
	        0 & * 
       \ebar \right\}, ~~~~~~
\mathcal{V}_1 = \left\{ \bbar {cc}
  	       0 & * \\
	       {*} & 0
       \ebar \right\}.
\edm
As in Section \ref{systemsoffinitetype} we can define a Lie algebra
of $so(2n,{\mathbb C})[z, z ^{-1}]$ by
\bdm
\gtildesig := \{ \sum_{i=\alpha}^{\beta} X_i z ^i ~|~ 
X_i \in \mathcal{V}_{i ~ \textup{(mod2)}} \}.
\edm
Observe that the Lie subalgebra $\mathcal{V}_0 \subset so(2n,{\mathbb C}) \subset 
\gtildesig$ has the Lie subalgebra decomposition:
\beqas
\mathcal{V}_0 &=& \mathcal{V}_0^U \oplus \mathcal{V}_0^L,\\
~~\\
\mathcal{V}_0^U = \left\{ \bbar {cc}
			{*} & 0 \\
			 0 & 0
	        \ebar \right\},  &&
\mathcal{V}_0^L = \left\{ \bbar {cc}
			 0 & 0 \\
			 0 & *
                 \ebar \right\},
\eeqas
Now let 
\bdm
\mathcal{P} := \{ X \in \gtildesig ~|~ X_0 \in \mathcal{V}_0^U,~
X_i = 0,~ \textup{for}~ i<0 \}.
\edm
It is clear that $\mathcal{P}$ is a Lie subalgebra of $\gtildesig$. 
The reason
we consider this subalgebra is that, comparing with (\ref{admissible}),
we see that admissible connections are precisely
those 1-forms whose coefficients are degree 1 elements of $\mathcal{P}$
and are real-valued for real values of $z$. The reality condition is 
captured in the prescription:
\beq \label{reality}
\overline{X(\bar{z})} = X(z).
\eeq
We denote by $\gtildesigr$ the subalgebra of $\gtildesig$ consisting of elements
which satisfy this condition.

 We now have the vector space decomposition:
\bdm
\gtildesig = \mathcal{P} \oplus \mathcal{N},
\edm
where $\mathcal{P}$, and 
$\mathcal{N} = \{X \in \gtildesig ~|~ X_0 \in \mathcal{V}_0^L, ~
 X_i = 0 ~\textup{for}~ i > 0 \}$ are
both Lie subalgebras consisting of polynomials in $z$ and $z^{-1}$
respectively.  This is exactly the type of decomposition 
to which the results of  Section \ref{systemsoffinitetype} apply. The
reality condition (\ref{reality}) is preserved by the Lax equation
(\ref{flowequation}) provided the polynomials involved have real
coefficients, so we 
can use Proposition \ref{loopgroupproposition} to produce admissible
connections 
\bdm
A:= \sum_{i=1}^n \pip \mathcal{V}_i(X) \dd t_i,
\edm
provided we can
find polynomials $\mathcal{V}_i(X) = p_i(X,z^{-1}) \in {\mathbb R}[X,z^{-1}]$
 which map $\gtildesig$
to $\gtildesig$ and are of top degree 1 in $z$.  In view of 
Remark \ref{conditionremark}  this can easily be
arranged using the fact that if $X$ is a skew-symmetric matrix then
so is $X^k$ for any odd positive integer, and we have the following
result:

\begin{proposition}  \label{solutionproposition}
The functions $\mathcal{V}_i: \gtildesig \to \gtildesig$ given by
\bdm
\mathcal{V}_i(X) = z^{2-2i} X^{2i-1},~i=1...n,
\edm
satisfy the condition (\ref{conditionA}) of Proposition \ref{loopgroupproposition}
for the loop algebra $\gtildesig$ defined above. This choice of $\mathcal{V}_i$ 
gives for each real initial condition
 $X_0(z) = \overline{X_0(\bar{z})} \in (\gtildesigr)^1_d$ a solution
$X(t)$ of (\ref{lax}) such that:
\begin{enumerate}
\item
 $A(X)$ is an admissible connection.
\item
  $X: {\mathbb R}^n \to \gtildesig$ is complete, and so, therefore is a
solution $F: {\mathbb R}^n \to SO(2n)$ obtained by integrating $A(X)$
evaluated at $z= z_0 \in {\mathbb R} _*$.
\end{enumerate}
\end{proposition}

\subsection{Remarks Concerning Proposition \ref{solutionproposition}}
\subsubsection{Completeness}
Flows of the form $dX = [X, V(X)]$, where $V(X)$ takes its values in
the same loop Lie algebra $\gtildesig$ as $X$, will always be complete if 
the underlying real Lie algebra is compact semi-simple (such as
 $so(n)$).  This is because in this case there is a positive definite
ad-invariant inner product $\leftb \cdot ~,~ \cdot \rightb$ on $\gtildesig$, where ad-invariant means
\bdm
\left< \, [A,B],C \,\right> = \leftb [C,A],B \rightb
\edm
 for any $A$, $B$ and $C$. Thus 
\beqas
\nabla \leftb X,X \rightb  &=& 2 \leftb \nabla X, X \rightb \\
&=& 2 \leftb [X, V(X)],X \rightb \\
& =&  2 \leftb [X,X], V(X) \rightb  \\
& =& 0,
\eeqas
 which is to say that $\leftb X,X \rightb $ is constant. 
Hence the solution $X$ is bounded, and
 therefore complete. 

\subsubsection{The Immersion Property} 
Every initial condition $X(z,0)$ in $(\gtildesig)_d$ will 
generate a family of complete maps $f_{z}: {\mathbb R} ^n \to S^{2n-1}$,
for $z \in {\mathbb R}$.  The map $f_{z}$ fails to be an immersion 
at points $t$
where the first $n$ entries of the $(n+1)$'th row of $F$ fail
 to be linearly independent 1-forms.
If we define $M(X)$ to be the $n \times n$ matrix whose $i$-th row
is the non-zero part of the $(n+1)$'th row of $X_1^{2i-1}$, namely 
\bdm
M(X) := \bbar {cccc}
(X_1)_{(n+1)1} & (X_1)_{(n+1)2} & ... & (X_1)_{(n+1)n} \\
(X_1^3)_{(n+1)1} & (X_1^3)_{(n+1)2} & ...  & (X_1)_{(n+1)n} \\
\vdots & \vdots & \vdots & \vdots \\
(X_1^{2n-1})_{(n+1)1} & (X_1^{2n-1})_{(n+1)2} & ... & (X_1^{2n-1})_{(n+1)n} 
\ebar,
\edm
then the condition for $f_{z}$ to be an immersion is 
\beq  \label{immersioncondition}
\det (M(X(z,t))) \neq 0.
\eeq
While most initial conditions $X(z,0)$ will ensure that this 
holds in a neighbourhood of $t=0$, in general one should expect to 
come across values of $t$ where the condition fails.

\subsubsection{Parallel Frames and Curved Flats}  \label{parallelremark}
As an alternative, one can define 
\beqas
\mathcal{P}_a & :=& \{ X \in \gtildesig ~|~ X_i = 0,~ i < 1 \}, \\
\mathcal{N}_a &:=& \{ X \in \gtildesig ~|~ X_i = 0, ~ i>0 \}.
\eeqas
These are again subalgebras, and so the same theory can be used
to produce admissible connections which have the additional property
that $\omega$ and $\eta$ are both zero, due to the fact that 
$\pip (... + Y_0 + Y_1 z) = Y_1 z$. 
Such connections are associated to frames which are 
parallel.  

On the other hand, if we were to set
\beqas
\mathcal{P}_a & :=& \{ X \in \gtildesig ~|~ X_i = 0,~ i < 0 \}, \\
\mathcal{N}_a &:=& \{ X \in \gtildesig ~|~ X_i = 0, ~ i \geq 0 \},
\eeqas
we would have the \emph{curved flats} which were studied in 
\cite{feruspedit1996II}.  These give you moving frames
 $F= [e_1,...,e_n, \xi_1,...\xi_n]$ which satisfy the flatness
conditions $\dd \omega + \omega \wedge \omega = 0$ and
$\dd \eta + \eta \wedge \eta =0$. Given a solution $F$ of this system,
one can parallelize the tangent and normal bundle by right multiplication
by a matrix $G$. One can show that $G$ is constant in $z$, but not in $t$,
and that $FG$ is actually the solution of the scheme for producing parallel
frames corresponding to the same initial condition $X(0)$.   
 
\subsubsection{Periodicity}
As stated in the introduction,
it is an interesting question to look for solutions, $f$, which 
 have $n$ linearly independent periods, corresponding to 
 flat $n$-tori in $S^{2n-1}$. For the case of parallel frames,
 because the frames, $F$,
are parallel for the immersion $f$, it follows
from the flatness of $f$ (parallel transport is independent of path)
that $f(x) = f(y)$ is equivalent to $F(x) = F(y)$, thus the periodicity
question for the immersion $f$ is precisely the periodicity question for the adapted frame $F$. Since parallel frames come from  simple $k$-symmetric
systems, we can use the results of Section \ref{properties} to say
that periodicity for the line bundle implies type two quasiperiodicity
for $F$.

 \textbf{Example 1} 
The simplest examples are of the form $X(0) = X_1 z $ which
give the same solutions whether we use the projections to 
$\mathcal{P}$ or to $\pa$, namely
the constant solutions $X(t) =  X(0),$
 because $[X, \pip \mathcal{V}_i (X)]$ are zero. 
For the case $n=2$, the most general such example is of 
the form
\bdm
X(z,t) = X(z,0) = \bbar {cccc}
        0 & 0 & x_1 & x_2 \\
	0 & 0 & y_1 & y_2 \\
	-x_1 & -y_1 & 0 & 0 \\
	-x_2 & -y_2 & 0 & 0
     \ebar z .
\edm
This generates the connection
\beqa
A&:=& X \dd t_1 + z^{-2} X^3 \dd t_2 \nonumber \\
&=&  \bbar {cccc}
        0 & 0 & x_1 & x_2 \\
	0 & 0 & y_1 & y_2 \\
	-x_1 & -y_1 & 0 & 0 \\
	-x_2 & -y_2 & 0 & 0
     \ebar z \, \dd t_1  +
\bbar {cc}
 0 & -C \\
 C^T & 0 
\ebar z \, \dd t_2, \label{simplecon}
\eeqa
where $C$ is the matrix
\bdm
\bbar {cc}
 x_1(x_1^2 +y_1^2) + x_2(x_1x_2 + y_1y_2)  & ~~   
     x_1(x_1x_2 + y_1y_2) + x_2(x_2^2 + y_2^2) \\
y_1(x_1^2 + y_1^2) + y_2(x_1 x_2 + y_1y_2) & ~~
      y_1(x_1x_2 + y_1y_2) +y_2(x_2^2 + y_2^2) 
\ebar.
\edm

The immersion condition (\ref{immersioncondition}) in this case is
\beqas
\left| \begin{array} {cc}
     x_1 & y_1 \\
    x_1(x_1^2 + y_1^2) + x_2(x_1x_2 + y_1y_2) & ~~
      y_1(x_1^2 + y_1^2) + y_2(x_1x_2 + y_1y_2)
  \end{array} \right|  \\
~~~~~ = (x_1x_2 + y_1y_2)(x_1y_2-y_1x_2)\\
~~~~~ \neq  0. ~~~~~~~~~~~~~~~~~~~~~~~~~~~~~~~~~~~~~
\eeqas
Any matrix $X$ of the above form gives rise to a complete
flat immersion $f: {\mathbb R}^2 \to S^3$, provided this condition holds. 

Since the second fundamental form is constant for such a connection, one
can perform a global change of coordinates so that it is diagonal and
constant everywhere, which means that the submanifold is a product of circles,
also known as a Clifford torus. This argument also applies
to  $n >2$.

Conversely, we can show that any Clifford 2-torus comes from the scheme
in the following way: a Clifford 2-torus has the parameterization 
\bdm
f(s_1,s_2) = [a \cos s_1, \, a \sin s_1, \, b \cos s_2, \, b \sin s_2],
\edm
where $a^2 + b^2 =1$.
We can choose the following adapted frame for $f$:
\bdm
F := \left(  \begin{array} {cccc}
  -\sin s_1 & 0 & a \cos s_1 & b \cos s_1 \\
 \cos s_1 & 0 & a \sin s_1 & b \sin s_1 \\
 0 & -\sin s_2 & b \cos s_2 & - a \cos s_2 \\
 0 & \cos s_2 & b \sin s_2 & -a \sin s_2
\end{array} \right),
\edm
and calculate
\beq  \label{cliffordcon}
F^{-1} \dd F = \left( \begin{array} {cccc}
  0 & 0 & a & b \\
  0 & 0 & 0 & 0 \\
  -a & 0 & 0 & 0 \\
  -b & 0 & 0 & 0
   \end{array} \right) \dd s_1
+
\left( \begin{array} {cccc}
0 & 0 & 0 & 0 \\
0 & 0 & b & -a \\
0 & -b & 0 & 0 \\
0 & a & 0 & 0 
\end{array} \right) \dd s_2.
\eeq

It is straightforward from here to solve the equation
\bdm
\bbar {c}
 s_1 \\
 s_2
\ebar 
 = B 
\bbar {c}
 t_1 \\
t_2
\ebar,
\edm
for the matrix $B$ so that after a linear change of coordinates the
connection (\ref{cliffordcon}) has the form of the right hand
side of (\ref{simplecon}).




\end{document}